\documentclass[10pt,a4paper,epsfig,final]{siamltex}
\usepackage{epsfig,longtable,subfigure}
\usepackage{amsmath,amssymb,float}
\floatstyle{boxed}
\newfloat{algorithm}{thp}{loa}
\floatname{algorithm}{Algorithm}
\def \p { \psi }

\def \f {\phi}

\title{Numerical Computations of Viscous, Incompressible 
Flow Problems Using a Two-Level Finite Element Method}
\author{Faisal Fairag\thanks{Mathematics Department, King Fahd 
University of Petroulem \& Minerals, Dhahran, 31261 Saudi Arabia 
         ({\tt ffairag@kfupm.edu.sa}).} }

\begin{document}

\maketitle

\begin{abstract}
We consider two-level finite element discretization methods for the 
stream function formulation of the Navier-Stokes equations.
The two-level method consists of solving a small nonlinear system on
the coarse mesh, then solving a linear system on the fine mesh. The
basic result states that the errors between the coarse and fine
meshes are related superlinearly.  This paper demonstrates that the
two-level method can be implemented to approximate efficiently
solutions to the Navier-Stokes equations.  Two fluid flow
calculations are considered to test problems which have a known
solution and the driven cavity problem.
Stream function contours are displayed showing the main features of
the flow.
\end{abstract}



\pagestyle{myheadings}
\thispagestyle{plain}

\section{Introduction}

The numerical treatment of nonlinear problems that arise in areas
such as fluid mechanics often requires solving large systems of
nonlinear equations.  Many methods have been proposed that attempt to
solve these systems efficiently; one such class of methods are
two-level methods that reduce the computational time.  The 
computational attractions of the methods are that they require the 
solution of only a small system of nonlinear equations on coarse mesh 
and one linear system of equations on fine mesh.  Apparently, the 
two-level method was proposed first in \cite{xu92,xu94a,xu94b} and 
used for semilinear elliptic problems. The method was implemented for 
the velocity-pressure formulation of the Navier-Stokes equations in 
\cite{layton93,layton95,layton96} and for the stream function 
formulation of the Navier-Stokes equations in \cite{fairag,ye,myphd}.

The Navier-Stokes equations may be solved using either the primitive
variable or stream function formulation. Here we use the stream
function formulation.  The attractions of the stream function
formulation are that the incompressibility constraint is
automatically satisfied, the pressure is not present in the weak
form, and there is only one scalar unknown to solve for. The standard
weak formulation of the stream function version first appeared in
1979 in \cite{gr79}. In this direction, Cayco and Nicolaides
\cite{cn86,cn89,cayco} studied a general analysis
of convergence for this standard weak formulation.

The goal of this paper is to demonstrate that the two-level method
can be implemented to approximate solutions for incompressible viscous
flow problems with high Reynolds number.

\section{Governing Equations}

Consider the Navier-Stokes equations describing the flow of
an incompressible fluid
% -----------------------
\begin{alignat}{3}
    - {\rm Re}^{-1}\bigtriangleup \vec{u} + (\vec{u} \cdot \bigtriangledown)\vec{u} + 
	\bigtriangledown p &= \vec{f} &\qquad \mbox{in} &\qquad \Omega, \label{nse1}\\
	\bigtriangledown \cdot \vec{u} &= 0 &\qquad \mbox{in} &\qquad \Omega, 
	\label{nse2} \\
	\vec{u} &= 0 &\qquad \mbox{on} &\qquad \partial \Omega, \label{nse3}
\end{alignat}
% ------------------------
where $\vec{u} = (u_1, u_2)$ and $p$ denotes the unknown velocity and
pressure field, respectively, in a bounded,simply connected polygonal
domain $\Omega \subseteq R^2$. $\vec{f}$ is a given body force and Re 
is the Reynolds number.

The introduction of a stream function $\psi(x,y)$ defined by
\begin{equation*}
u_1 = - \frac{\partial \psi}{\partial y}, \qquad u_2 =
\frac{\partial \psi}{\partial x} ,
\end{equation*}
means that the continuity equation (\ref{nse2})
is satisfied identically.
The pressure may then be eliminated from (\ref{nse1}) to give
%-------------------------
\begin{alignat}{3}
\mbox{Re}^{-1}\triangle^2 \psi - \psi_y \triangle \psi_x +
\psi_x \triangle \psi_y & = & ~\vec{{\rm curl}} \vec{f}
 \mbox{ in }  \Omega , \label{sf1} \\
\psi & = & 0  \mbox{ on }  \partial \Omega , \label{sf2} \\
\frac{\partial \psi}{\partial \hat{n}} & = & 0  \mbox{ on } 
\partial \Omega , \label{sf3}
\end{alignat}
%-------------------------
where $\hat{n}$ represents the outward unit normal to $\Omega$.
In order to write (\ref{sf1})-(\ref{sf3}) in a variational form, we
define the
Sobolev spaces
%-------------------------
\begin{alignat}{4}
H^1(\Omega)  & =  \{v: v\in L^2(\Omega), Dv \in L^2(\Omega)\}, \\
H^1_0(\Omega)& =  \{v: v\in H^1(\Omega), v = 0 \mbox{ on } \partial \Omega\},\\
H^2(\Omega) & =  \{v: v\in L^2(\Omega), Dv\in L^2(\Omega), D^2v \in L^2(\Omega)\}, \\
H^2_0(\Omega) & =  \{v\in H^2(\Omega): v = \frac{\partial v}{\partial n} = 0,
\mbox{ on } \partial \Omega\},
\end{alignat}
%-------------------------
where $L^2(\Omega)$ is the space of square integrable functions on
$\Omega$ and $D$ represents differentiation with respect to $x$ or 
$y$. For each $\phi\in H^1(\Omega)$, define $\vec{{\rm curl}} \; \phi 
= \left(\begin{array}{c} \phi_y \\
-\phi_x\end{array} \right)$. The standard weak form of equations
(\ref{sf1})-(\ref{sf3}) is:
%-------------------------
\begin{equation}   \label{cont_stream_form}
\begin{split}
  \mbox{Find } \psi \in H_{0}^{2}(\Omega) \mbox{ such that, for all }
  \phi \in 
  H_{0}^{2}(\Omega), \\
  a(\psi,\phi)+b(\psi;\psi,\phi) = l(\phi), 
\end{split}
\end{equation}
%-------------------------
where
\begin{equation*}
\begin{split}
    a(\p,\f) &= Re^{-1}\int_{\Omega}\bigtriangleup
    \p \cdot \bigtriangleup \f ~d\Omega  ,      \\
    b(\xi;\p,\f) &= \int_{\Omega} \bigtriangleup
    \xi(\p_{y} \f_{x}-\p_{x}\f_{y})~d\Omega  ,  \\
    l(\f) &= (\vec{f}, \vec{\mbox{curl}\f}) 
	= \int_{\Omega} \vec{f} \cdot 
    \vec{\mbox{curl}} \f ~d\Omega .                                                                 
\end{split}
\end{equation*}

\section{Finite Element Discretization}

For the standard finite element discretization of
(\ref{cont_stream_form}) we choose
subspace $X^H \subset H^2_0(\Omega)$. We then 

\begin{equation} \label{disc_stream_form}
\begin{split}
  \mbox{seek } \psi^H \in X^H \mbox{ such that, for all } \phi^H \in 
  X^H, \\
  a(\psi^H,\phi^H)+b(\psi^H;\psi^H,\phi^H) = l(\phi^H). 
\end{split}
\end{equation}

One can prove existence and uniqueness for the solution of the
discrete problem (\ref{disc_stream_form}) ; see \cite{cn86} for low
Reynolds numbers and
see \cite{fairag,myphd} for high Reynolds numbers.

Once the finite element spaces are prescribed, the discrete
problem (\ref{disc_stream_form}) reduces to solving a system of
nonlinear algebraic
equations which has a Jacobian which is large, sparse and bounded.
Various iterative methods can be used to solve the nonlinear problem 
(\ref{disc_stream_form}). For example, a standard approach is to use
Newton's method to
linearize (\ref{disc_stream_form}) for a fixed mesh of size $h$.
Since only one fixed
mesh spacing is used, we will refer to this approach as a one-level
approach. In the next section we will describe a two-level method for
solving the discrete problem (\ref{disc_stream_form}).

\section{Two-level Method}

We consider the approximate solution of (\ref{sf1}) by a two-level,
finite-element procedure.  Let $X^h, X^H \subset H^2_0(\Omega)$
denote two conforming finite-element meshes with $H \gg h$.  The 
method we consider computes an approximate solution $\psi^h$ in the
finite-element space $X^h$ by solving one linear system for the
degrees of freedom on $X^h$.  This particular linear problem requires 
the construction of a finite-element space $X^H$ upon a very coarse
mesh of width '$H \gg h$', and then the solution of a much smaller 
system of nonlinear equations for an approximation in $X^H$.  The 
solution procedure is then given in Algorithm \ref{two_level_alg}.
%------------test test test test --------------------
%\vspace{4mm}
%\noindent \hspace*{7mm} {\sc \bf Algorithm}.

%\vspace{4mm}
\begin{algorithm*} \caption{The Two-Level Algorithm}  \label{two_level_alg}
\noindent \hspace*{7mm} {\bf Step 1.} Solve the nonlinear system on coarse mesh for
$\psi^H \in X^H$:
\begin{equation}
a\left(\psi^H, \phi^H\right) + b\left(\psi^H, \psi^H, \phi^H\right)
= \left(\vec{f}, \vec{{\rm curl}} \;
\phi^H\right), \hspace*{7mm} \mbox{ for all } ~ \phi^H\in X^H.
\end{equation}
\hspace*{7mm} {\bf Step 2.} Solve the linear system on fine mesh for $\psi^h \in X^h$:
\begin{equation}
a\left(\psi^h, \phi^h\right) + b\left(\psi^H, \psi^h, \phi^h\right)
= \left(\vec{f}, \vec{{\rm curl}} \;
\phi^h\right), \hspace*{7mm} \mbox{ for all } ~ \phi^h\in X^h.
\end{equation}
\end{algorithm*}
%------------test test test test --------------------
The inclusion $X^H \subset H^2_0(\Omega)$ requires the
use of finite-element functions that are continuously
differentiable over $\Omega$. We shall give some examples
of finite-element spaces for the stream function formulation~
(see ~\cite{ciarlet,cayco}).  We will impose boundary conditions by
setting all the degrees of freedom at the boundary nodes to be 
zero and the normal derivative equal to zero at all vertices
and nodes on the boundary.

\vspace{4mm}

\noindent {\sc \bf Argyis Triangle}. The functions are quintic
polynomials
within each triangle and the 21 degrees of freedom are chosen
to be the function value, the first and second derivatives at the
vertices, and the normal derivative at the midsides.

\vspace{4mm}

\noindent {\sc \bf Clough-Tocher Triangle}. Here we subdivide each
triangle into
three
triangles by
joining the vertices to the centroid.  In each of the smaller
triangles, the functions are cubic polynomials.  There are then 30
degrees of freedom needed to determine the three different cubic
polynomials associated with the three triangles.  Eighteen of these
are
used to ensure that, within the big triangle, the functions are
continuously differentiable.  The remaining 12 degrees of freedom
are
chosen to be the function values and the first derivatives at the
vertices and the normal derivative at the midsides.

\vspace{4mm}

\noindent {\sc \bf Bogner-Fox-Schmit Rectangle}. The functions are
bicubic
polynomials within each rectangle. The degrees of freedom are
chosen to
be the function value, the first derivatives, and the mixed second
derivative at the vertices.  We set the function and the normal
derivative values equal to zero at all vertices on the boundary.

\vspace{4mm}

\noindent {\sc \bf Bicubic Spline Rectangle}. The functions are the
product
of
cubic splines.  These functions are bicubic polynomials within each
rectangle, are twice continuously differentiable over $\Omega$, and
their degrees of freedom are the function values at the nodes (plus some
additional ones on the boundary).

\vspace{4mm}

The question which automatically arises is how to choose $X^H$ and
$X^h$ so that we obtain optimal  accuracy.  This question was addressed 
from a theoretical standpoint in \cite{fairag} and the results are
summarized below. In \cite{fairag}, it was proven that the algorithm 
produces an approximate solution which satisfies the error bound
\begin{equation} \label{error}
\left|\psi-\psi^h\right|_2 \leq C\left\{
\inf_{w^h \in X^h} \left|\psi - w^h\right|_2 + \left|
\ln h\right|^{1/2}\cdot \left|\psi - \psi^H\right|_1\right\}.
\end{equation}
As an example, consider the case of the Clough-Tocher triangle.
For this
element we have the following inequalities:
\begin{eqnarray*}
\left|\psi -\psi^h\right|_j \leq Ch^{4-j} \hspace*{7mm} (j=0, 1,2),
\\
\left|\psi -\psi^H\right|_j \leq CH^{4-j} \hspace*{7mm} (j=0, 1,2).
\end{eqnarray*}
Thus, if we seek an approximate solution $\psi^h$ with the same
asymptotic accuracy as $\psi^h$ in $|\cdot|_2$, the above error
bound
shows that the superlinear scaling, between coarse and fine meshes,
\begin{equation}
h = O\left(H^{3/2}|\ln H|^{1/4}\right) ,
\end{equation}
suffices.  Analogous scalings between coarse and fine meshes can be
calculated from (\ref{error}) by balancing error terms on the
right-hand side of (\ref{error}) in the same way.  For each of the 
elements described above, we give, in Table~\ref{table1}, the scaling
between coarse and fine meshes.
%----------------test-----------
%----------------test-----------
\begin{table}[h]
\begin{tabular}{|c|c|c|c|} 
\multicolumn{4}{c}{} \\ 
\hline 
Element    & $\mid\psi-\psi^{H}\mid_{2}$ &
$\mid\psi-\psi^{H}\mid_{1}$ &
                     Scaling \\ \hline 
Argyris triangle               & $ H^{4} $ & $ H^{5} $ &
                     $ h \mid \ln h \mid^{-1/4} = O(H^{5/2}) $ \\
Clough-Tocher triangle         & $ H^{2} $ & $ H^{3} $ &
                     $ h \mid \ln h \mid^{-1/4} = O(H^{3/2}) $ \\
Bogner-Fox-Schmit rectangle    & $ H^{2} $ & $ H^{3} $ &
                     $ h \mid \ln h \mid^{-1/4} = O(H^{3/2}) $ \\
Bicubic spline rectangle       & $ H^{2} $ & $ H^{3} $ &
                     $ h \mid \ln h \mid^{-1/4} = O(H^{3/2}) $ 
                     \\\hline
\end{tabular}
\caption{Scaling between coarse and fine meshes}\label{table1}
\end{table}

\section{Numerical Examples}

In this section we describe some numerical results obtained by
implementing the two-level algorithm for which we have an exact
solution and the second is the well-known driven cavity problem.
We chose the later problem because there are numerous results in
the
literature with which to compare.

\vspace{4mm}

For both examples, the region $\Omega$ is the unit
square $\{ 0 < x < 1, \; 0 < y < 1\}$ and for the finite element
discretization we use the Bogner-Fox-Schmit elements.
In order to compare the efficiency of the proposed method all
linear and nonlinear systems were solved in the same way.  All
nonlinear problems were solved by Newton's method until the norm of
the difference in successive iterates and the norm of residual
were within a fixed tolerance. In each Newton's iteration, we need
to solve a linear system. The resulting linear system is
non-symmetric whose symmetric part is positive definite.
Moreover, the resulting matrix is sparse matrix.
We choose the Bi-Conjugate Gradient Stabilized method (BICGSTAB)
which requires two matrix-vector products and four inner products
in each iteration. BICGSTAB is given and discussed in~\cite{templates}.
When solving the linearized problem with a mesh spacing $h$ we need the
solution $u^H$ generated on a mesh with spacing $H$.  To do this we
interpolate the solution $u^H$ onto the grid with spacing $h$.

\vspace{4mm}

\noindent {\bf Example 1}

We consider as a test example the 2D Navier-Stokes equations (\ref{nse1})-
(\ref{nse3}) on the unit square $ \Omega = (0,1)^2 $ where we define the 
right hand side by $ f := -Re^{-1} \bigtriangleup \vec{u} +(\vec{u} 
\cdot {\bf
\bigtriangledown})
  \vec{u}+ \bigtriangledown  p $ with the following prescribed exact
  solution
\begin{eqnarray*}
u & = & \left(\begin{array}{c}
\psi_y \\ -\psi_x\end{array} \right) \hspace*{7mm} \mbox{with } \;
\psi(x,y) = x^2(x-1)^2 y^2(y-1)^2 , \\
p & = & x^3 + y^3 - 0.5.
\end{eqnarray*}
For this test problem, all requirements of the theory concerning 
the geometry of the domain and the smoothness of the data are
satisfied.
Moreover, the stream function $\psi(x,y)$ satisfies the boundary 
conditions of the stream function equation of the Navier-Stokes
equations.

\vspace{4mm}

\begin{figure}[H]
  \centerline{
  \psfig{figure=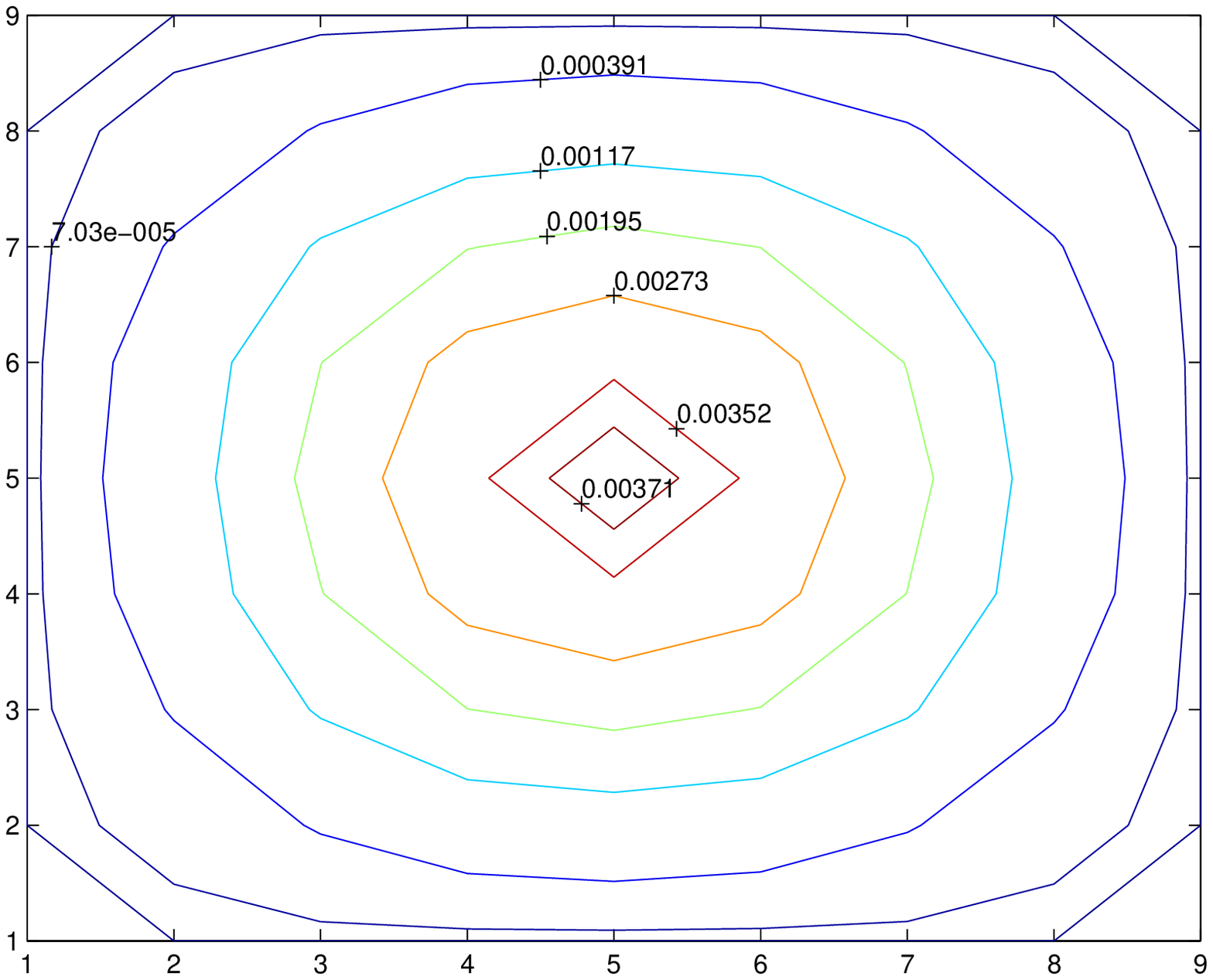,width=2.5in,height=2.5in}
  \psfig{figure=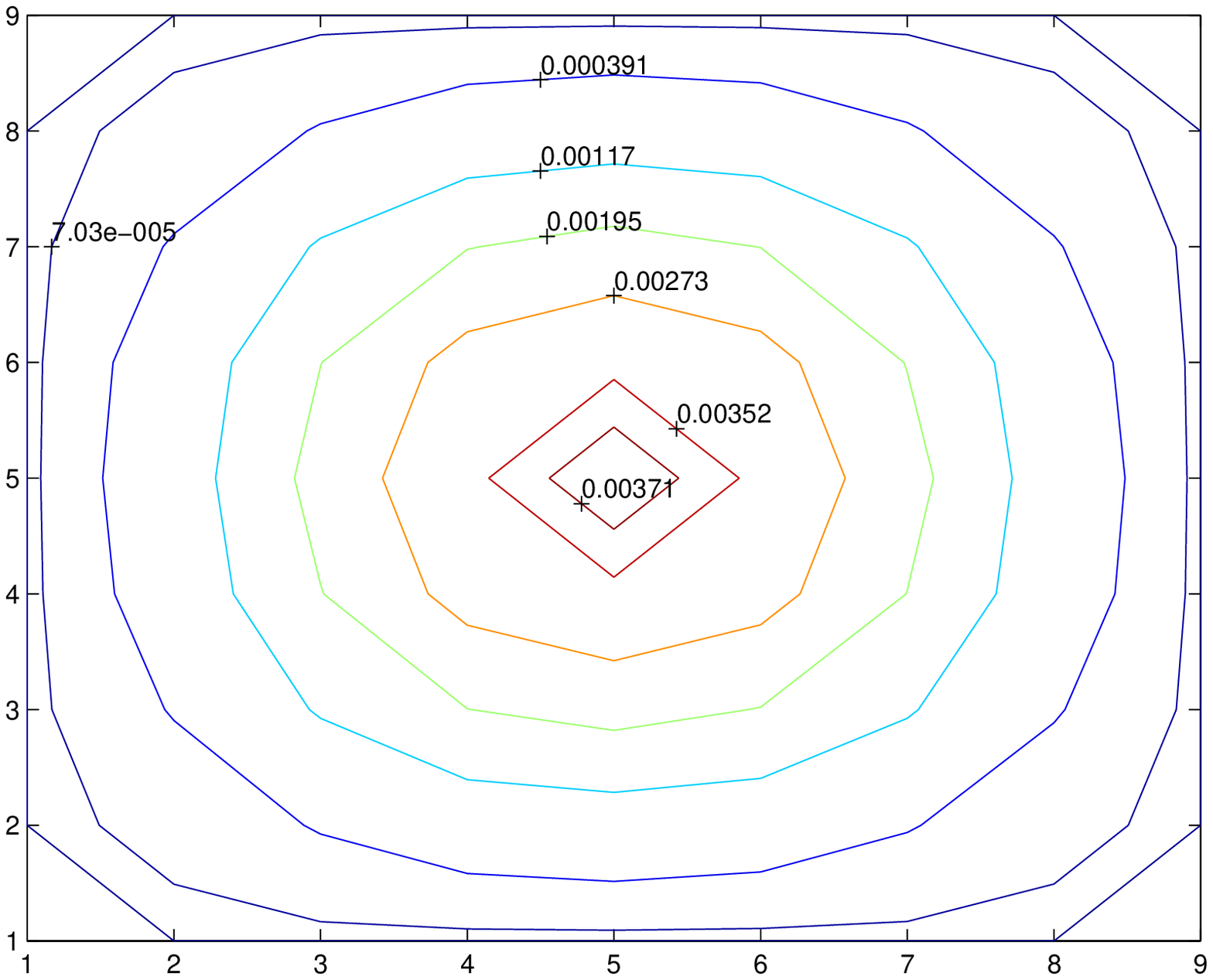,width=2.5in,height=2.5in}  }
  \centerline{
  \psfig{figure=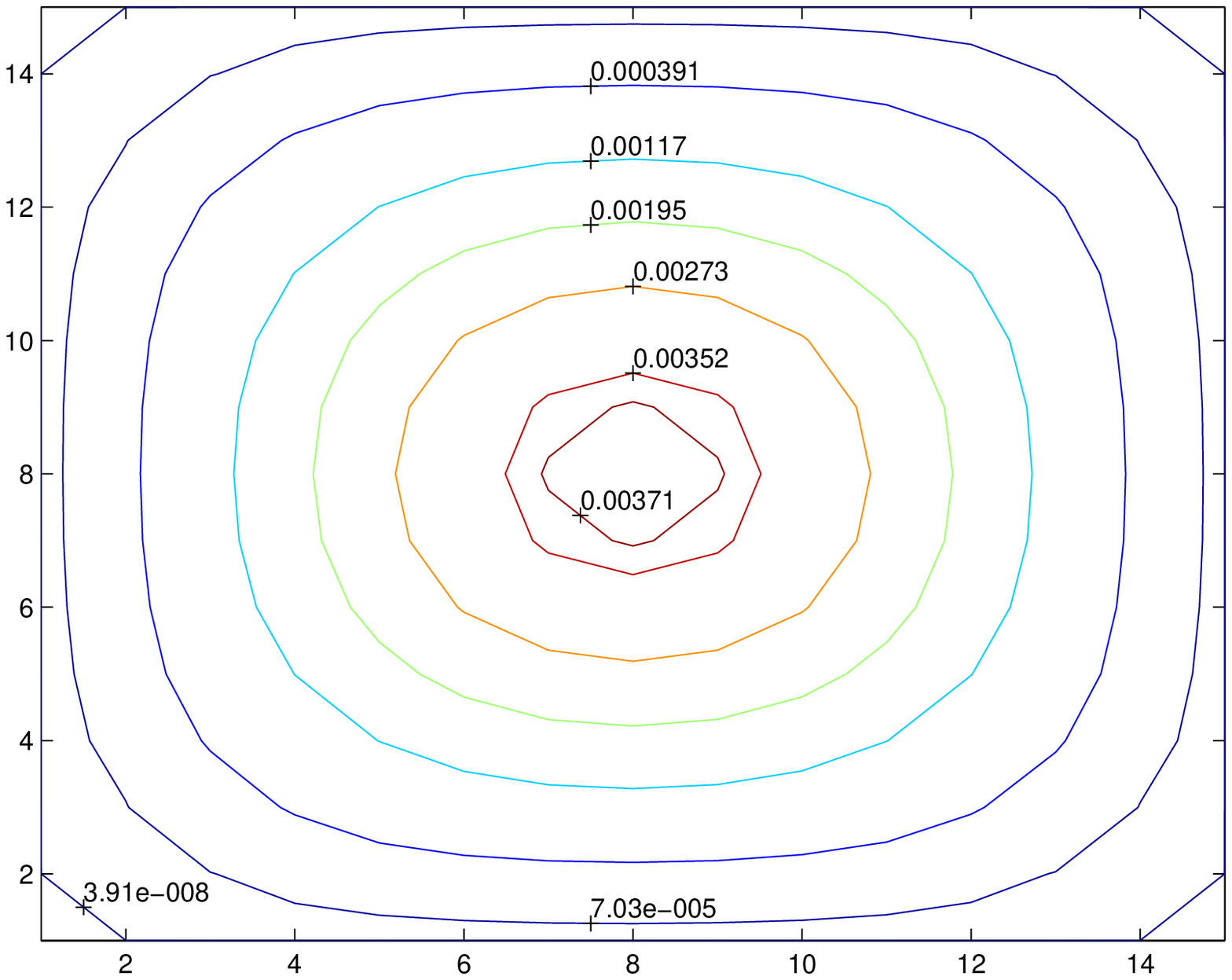,width=2.5in,height=2.5in}
  \psfig{figure=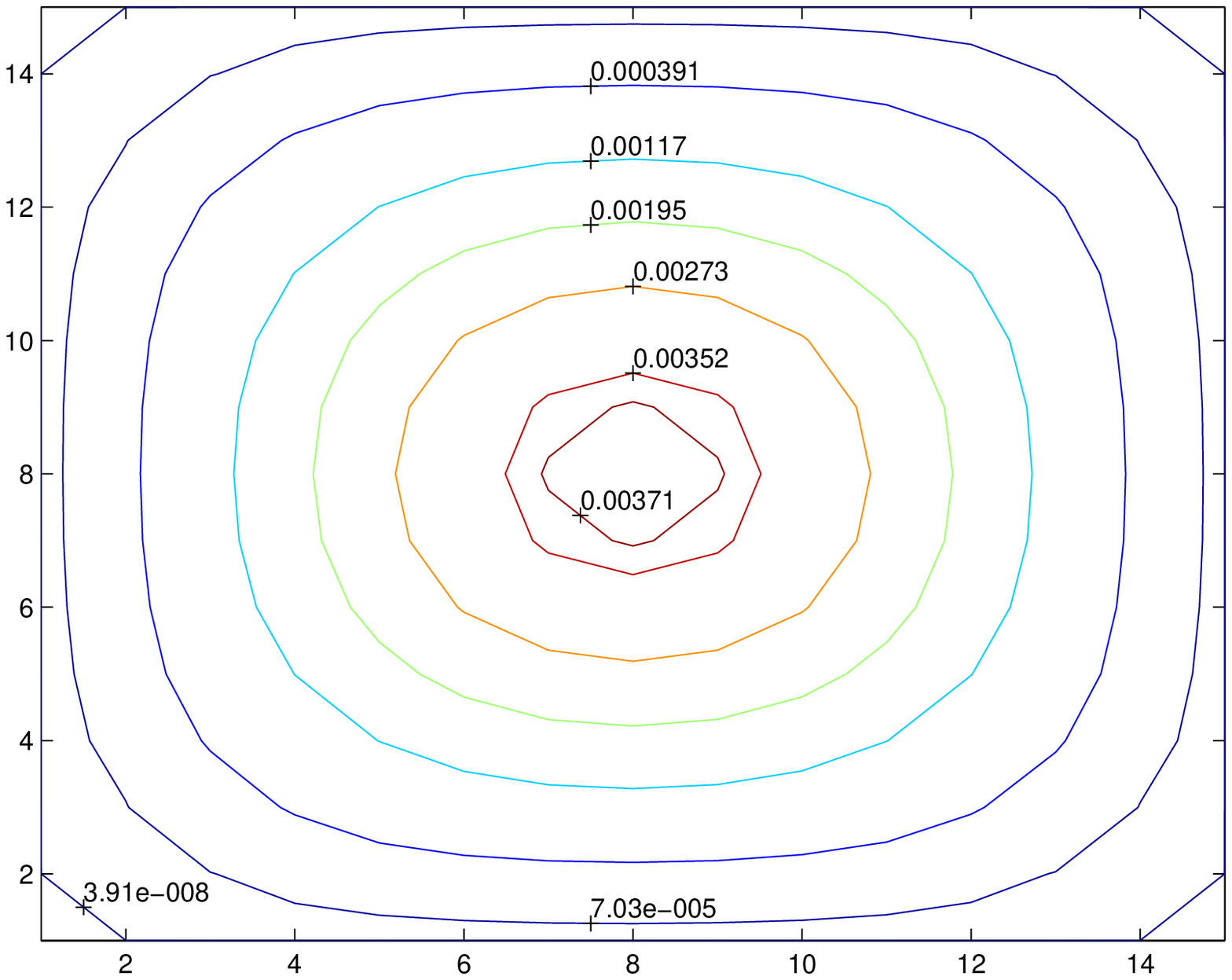,width=2.5in,height=2.5in}  }
  \centerline{
  \psfig{figure=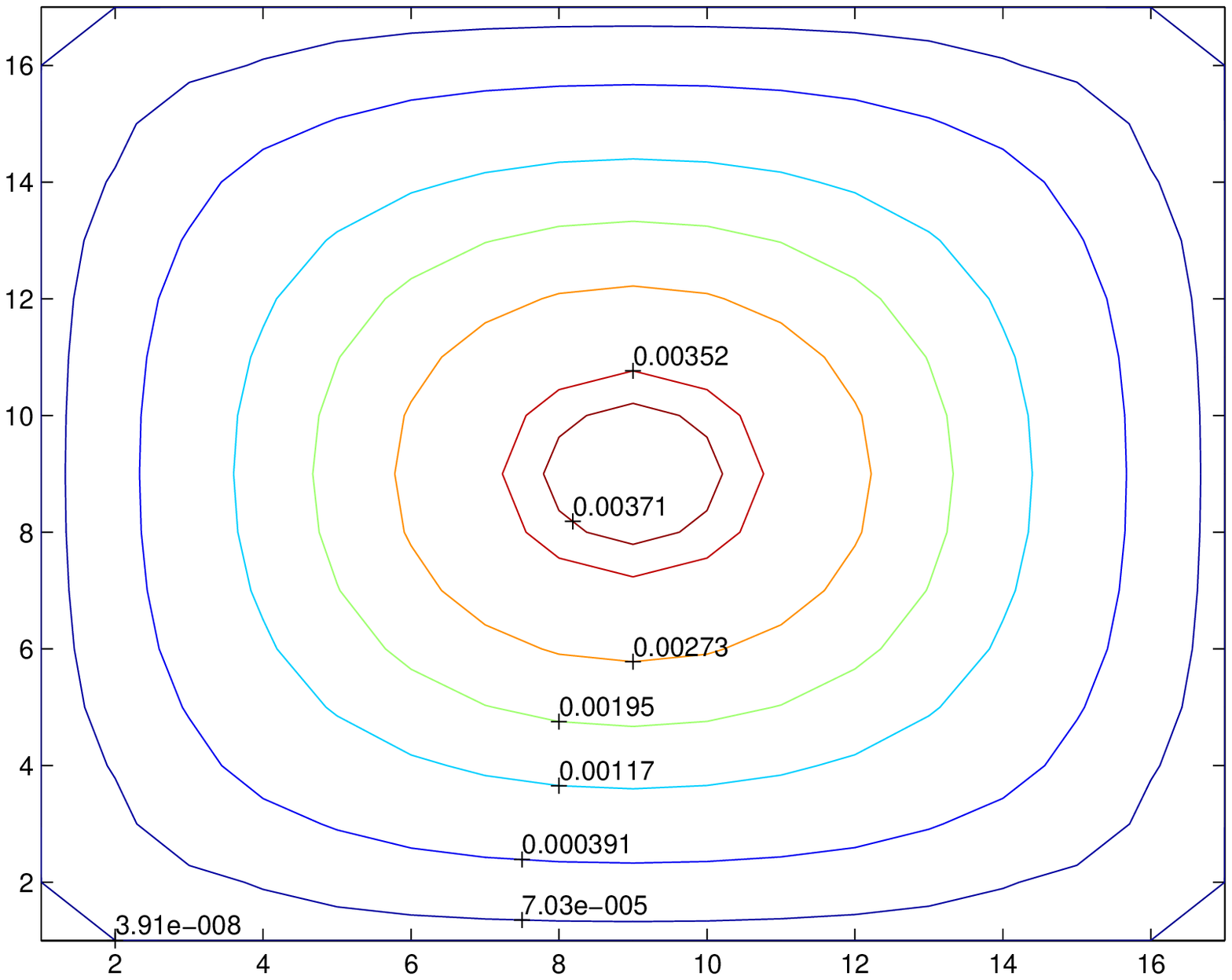,width=2.5in,height=2.5in}
  \psfig{figure=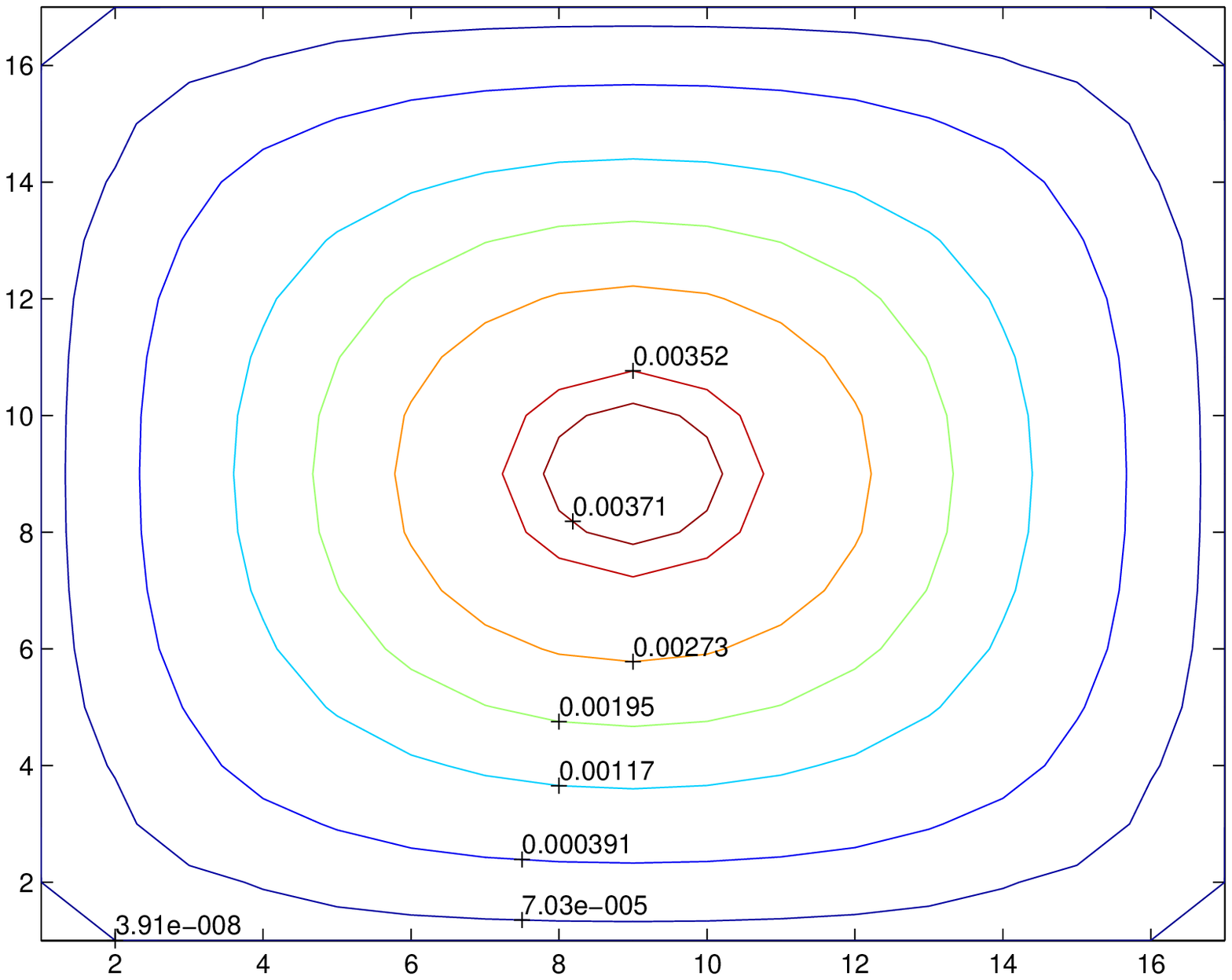,width=2.5in,height=2.5in}  }
  \caption{Left-side: streamlines for $h = \frac{1}{8},
           \frac{1}{14}, \frac{1}{16}$ with Re = 10 using 
		   the one-level method. Right-side: streamlines for 
		   $(H,h) = \left(\frac{1}{4}, \frac{1}{8}\right), 
		   \left(\frac{1}{7}, \frac{1}{14}\right), 
		   \left(\frac{1}{8}, \frac{1}{16}\right)$ 
		   with Re = 10 using the two-level method. }\label{fig1}
\end{figure}

Our goal in this test is to validate the code and the properties
and merits of the two-level method as compared with the one-level
method.  In all numerical calculations in this example we have
used the Bogner-Fox-Schmit elements with $\mbox{Re} = 10$ and
$\mbox{tol}
= 10^{-3}$.  We pick three values of $h$.  They are $\frac{1}{8},
\frac{1}{14}$ and $\frac{1}{16}$. The cpu-time, number of Newton's
iterations, number of Bicgstab iterations, the $L^2$-error and
$H^2$-error of the stream function $\psi$ for the one-level method
for
different values of $h$ are tabulated in Table~\ref{table2}.  
Table~\ref{table3} shows
cup-time, number of Newton's iterations and the number of
Bicgstab's
iterations for each linear solver for the two-level method.  
Figure~\ref{fig1}
shows, for fixed $\mbox{Re} = 10$ using the one-level and two-level
method, the streamlines for $h = \frac{1}{8}; \frac{1}{14};
\frac{1}{16}$.

\noindent {\bf Remarks}
\begin{enumerate}
\item From Tables~\ref{table2} and~\ref{table3}, the cpu-time for the
two-level
method is much smaller than the corresponding cpu-time for the
one-level
method.  For $h = \frac{1}{8}, \frac{1}{14}$, we save
about 57\%.  For example, in $h = \frac{1}{16}$ we
save about 73\%.  For example, in $h = \frac{1}{16}$ we
need to solve a
nonlinear system of 1156 equations which requires solving three
linear
systems of equations of order 1156.  The corresponding two-level
method
requires solving a nonlinear system of 324 equations and a linear
system of 1156 equations.  We anticipate the savings to increase as
the
mesh is further refined.
\item From Figures~\ref{fig1} , both columns are exactly the same,
which
means that the two-level method produces a solution with the same
quality as the one-level method.
\item From Tables~\ref{table2} and \ref{table3}, both $H^1$-error 
and $H^2$-error are of
the same order, which means that the velocity field is of the same
error and quality in both methods since $u = \psi_y$ and $v
=-\psi_x$.
\end{enumerate}
% --------------- Table 2 ----------------------
\begin{longtable}[h]{|c|c|c|c|c|c|c|}
\multicolumn{7}{c}{ } \\  
\hline  
$ h $   &  cpu time  & ni  & nb &
$\parallel \psi-\psi^{h}\parallel_{0,h}$ & $\parallel
\psi-\psi^{h}\parallel_{1,h}$ & $\parallel
\psi-\psi^{h}\parallel_{2,h}$  \\
\hline \hline 
$ \frac{1}{8} $ & 69.74 & 3 & 26,30,24 & 3.94$\times10^{-5}$ &
4.23$\times10^{-2} $ & 1.41$\times10^{-1}$ \\ \hline
$ \frac{1}{14} $ & 294.35 & 3 & 50,56,55 & 2.41$\times10^{-5}$ &
1.74$\times10^{-2}$ & 8.87$\times10^{-2}$ \\ \hline
$ \frac{1}{16} $ & 992.23 & 3 & 63,75,67 & 3.55$\times10^{-5}$ &
1.41$\times10^{-2}$ & 8.08$\times10^{-2}$ \\ \hline
\caption{One level method. where ni = number of Newton's iteration,
nb = number of Bicgstab iterations}\label{table2}
\end{longtable}
% --------------- Table 3 ----------------------
\begin{longtable}[h]{|c|c|c|c|c|c|c|}
\multicolumn{7}{c}{ } \\
\hline  
$ H,h $ & cpu time  & noc  & nof
& $\parallel \psi-\psi^{h}\parallel_{0,h}$ & $\parallel
\psi-\psi^{h}\parallel_{1,h}$ & $\parallel
\psi-\psi^{h}\parallel_{2,h}$  \\
\hline \hline 
$ \frac{1}{4} , \frac{1}{8} $ & 29.91 & 10,18,15 & 50 &
1.40$\times10^{-3}$ & 4.09$\times10^{-2} $ & 1.42$\times10^{-1}$ \\
\hline
$ \frac{1}{7} , \frac{1}{14} $ & 128.28 & 22,24,21 & 118 &
2.31$\times10^{-4}$ & 1.73$\times10^{-2}$ & 8.70$\times10^{-2}$ \\
\hline
$ \frac{1}{8} , \frac{1}{16} $ & 267.33 & 26,30,24 & 169 &
1.70$\times10^{-4}$ & 1.41$\times10^{-2}$ & 7.94$\times10^{-2}$ \\
\hline
\caption{Two level method.  where noc = number of Bicgstab
iterations in coarse, nof = number of Bicgstab iterations in fine}
\label{table3}
\end{longtable}
% ---------------Table 3 -------------------
\noindent {\bf  Example 2}

In this example, we consider the problem which was described
in Example 1.  The exact solution is very smooth and does not
depend on the Reynolds number.  The point of these tests is 
to increase the Reynolds number with $h$ fixed and test the 
robustness of the method.
Our goal in this test is to determine the validation of 
the code and the norm behavior when {\em Re} is varying.
Schieweck~\cite{schieweck} tested his code with this problem.
He used two nonconforming finite element approximations of
upwind type for the velocity-pressure formulation.

The numerical computations of this example were obtained using a
Sun Ultra 2 with 2 200\,Mhz ultrasparc processor running 
Solaris 2.5.1. In all numerical calculations, we used the 
Bogner-Fox-Schmit rectangles.
The streamlines for Re = 100, 1000 and 2000 were obtained with
$17 \times 17$ grid points on the coarse mesh and $33 \times 33$ grid
on the fine mesh.  Hence, a mesh of 256 elements and a mesh
of 1024 elements were used in this test for the case of
Bgner-Fox-Schmit
rectangles.

\vspace{4mm}

\begin{figure}[H]
\centering
\mbox{\subfigure[Exact]{\epsfig{figure=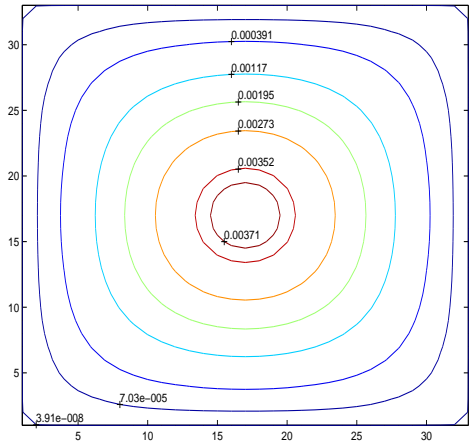,
      width=2.5in,height=2.5in}} \quad
      \subfigure[Re=100]{\epsfig{figure=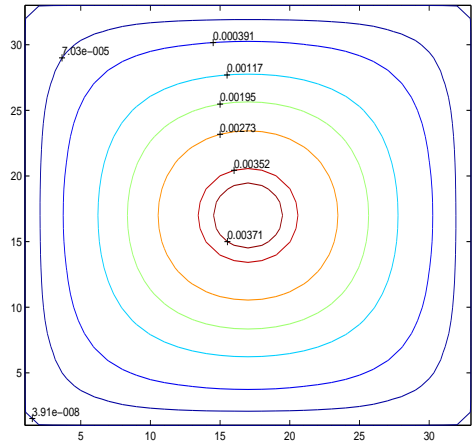,
	  width=2.5in,height=2.5in}} }
\mbox{\subfigure[Re=1000]{\epsfig{figure=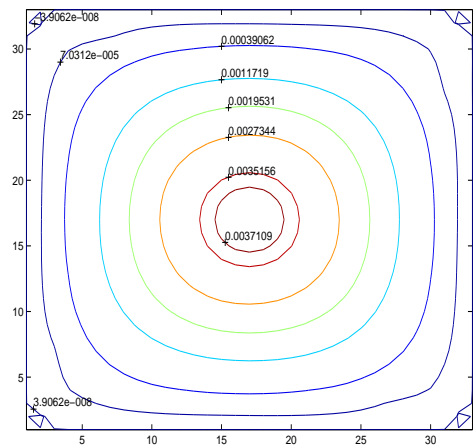,
	  width=2.5in,height=2.5in}} \quad
      \subfigure[Re=2000]{\epsfig{figure=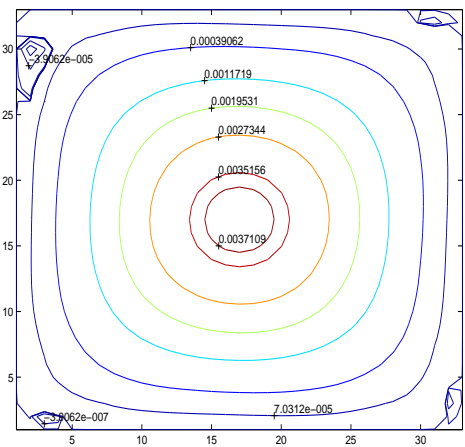,
	  width=2.5in,height=2.5in}}}
\caption{Streamlines for $ H= \frac{1}{16},h= \frac{1}{32} $ 
using Bogner-Fox-Schmit element. }\label{fig2}
\end{figure}

Table~\ref{table4} represents the $L^2$-error, $H^1$-error and
$H^2$-error.
The
streamlines are plotted in Figure~\ref{fig2} for Re = 100, 1000,
2000.
Figure
~\ref{fig2} shows that an  increase in the Reynolds number will
affect the
streamlines and increase the number of corner contours.
\vspace{4mm} \\
\noindent {\bf Example 3}
\vspace{4mm} \\
Cavity flows have been a subject of study for some time.
These flows have been widely used as test cases for validating
incompressible fluid dynamics algorithm. Corner
singularities for two-dimensional fluid flows are very
important since most examples of physical
interest have corners.  For example, singularities of most elliptic
problems develop when the boundary contour is not smooth.
In this example, we consider the driven flow in a rectangular
cavity when the top surface moves with a constant
velocity along its length.  The upper corners where the
moving surface meets the stationary walls are singular points of
the flow at which the horizontal velocity is multi-valued.
The lower corners are also weakly singular points.

\begin{figure}[H]
\centering
\mbox{\subfigure[Re=1]{\epsfig{figure=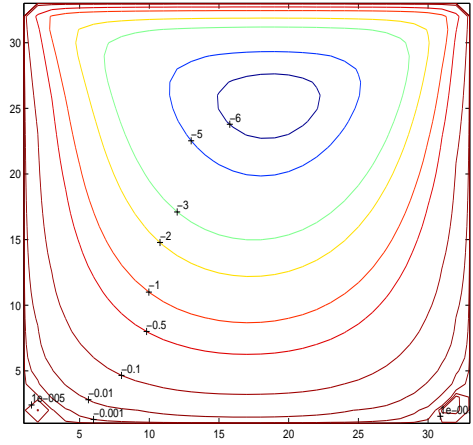,
	  width=2.5in,height=2.5in}} \quad
      \subfigure[Re=10]{\epsfig{figure=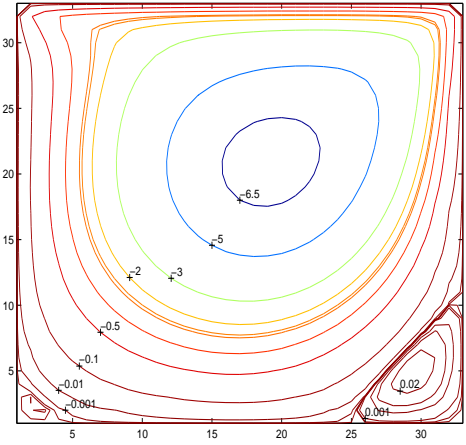,
	  width=2.5in,height=2.5in}}}
\mbox{\subfigure[Re=50]{\epsfig{figure=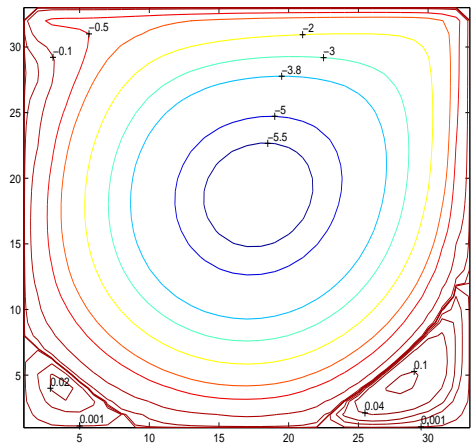,
	  width=2.5in,height=2.5in}} \quad
      \subfigure[Re=100]{\epsfig{figure=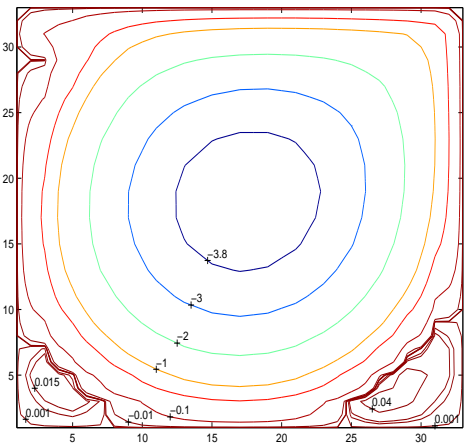,
	  width=2.5in,height=2.5in}}}	  
  \caption{Cavity Problem : Streamlines for $ H=\frac{1}{16}, h=
  \frac{1}{32} 
           $ with different values of Re numbers using two level 
		   method}\label{fig3}
\end{figure}

\begin{figure}[H]
\centerline{
\epsfig{figure=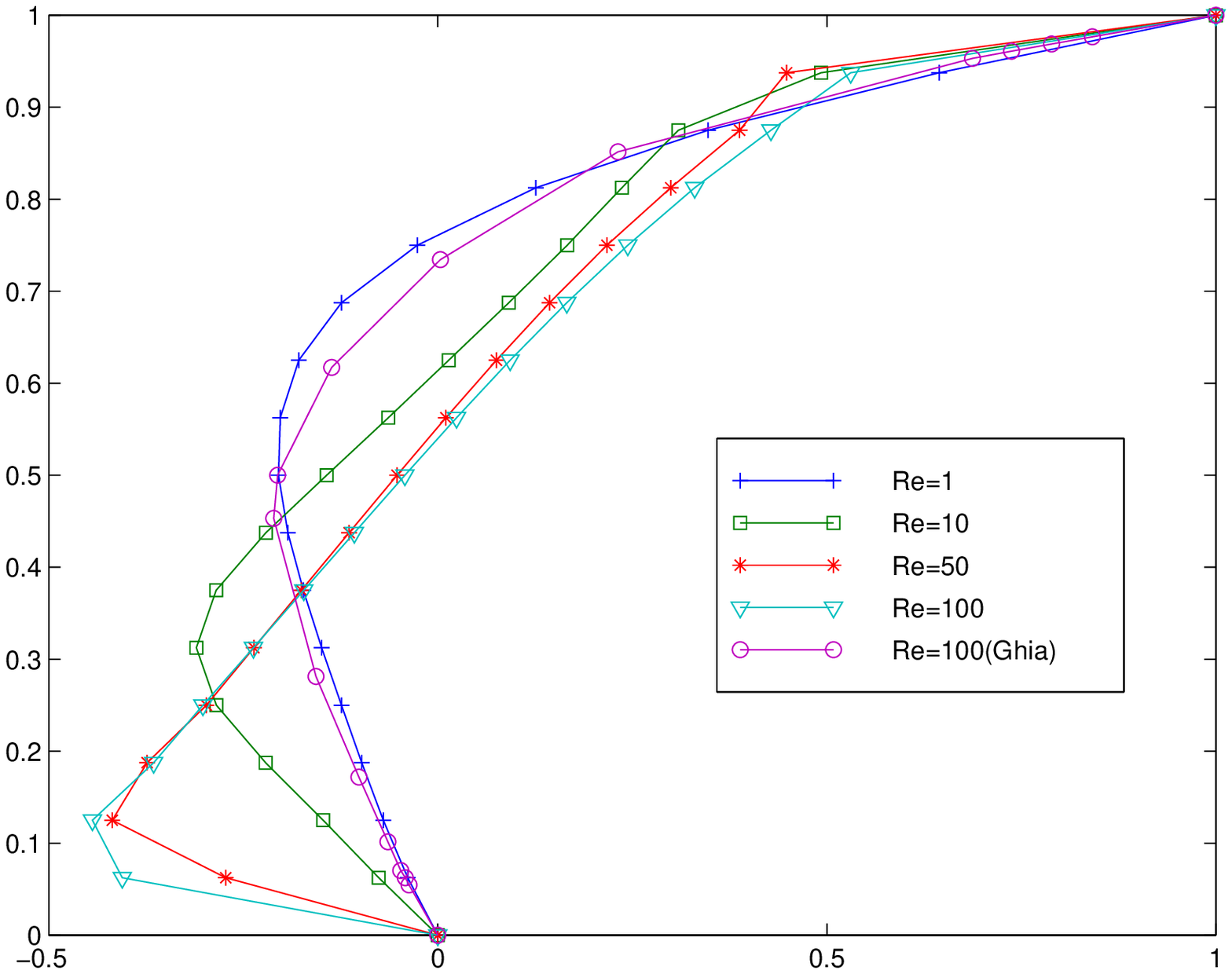,
	  width=2.5in,height=2.5in} 
\epsfig{figure=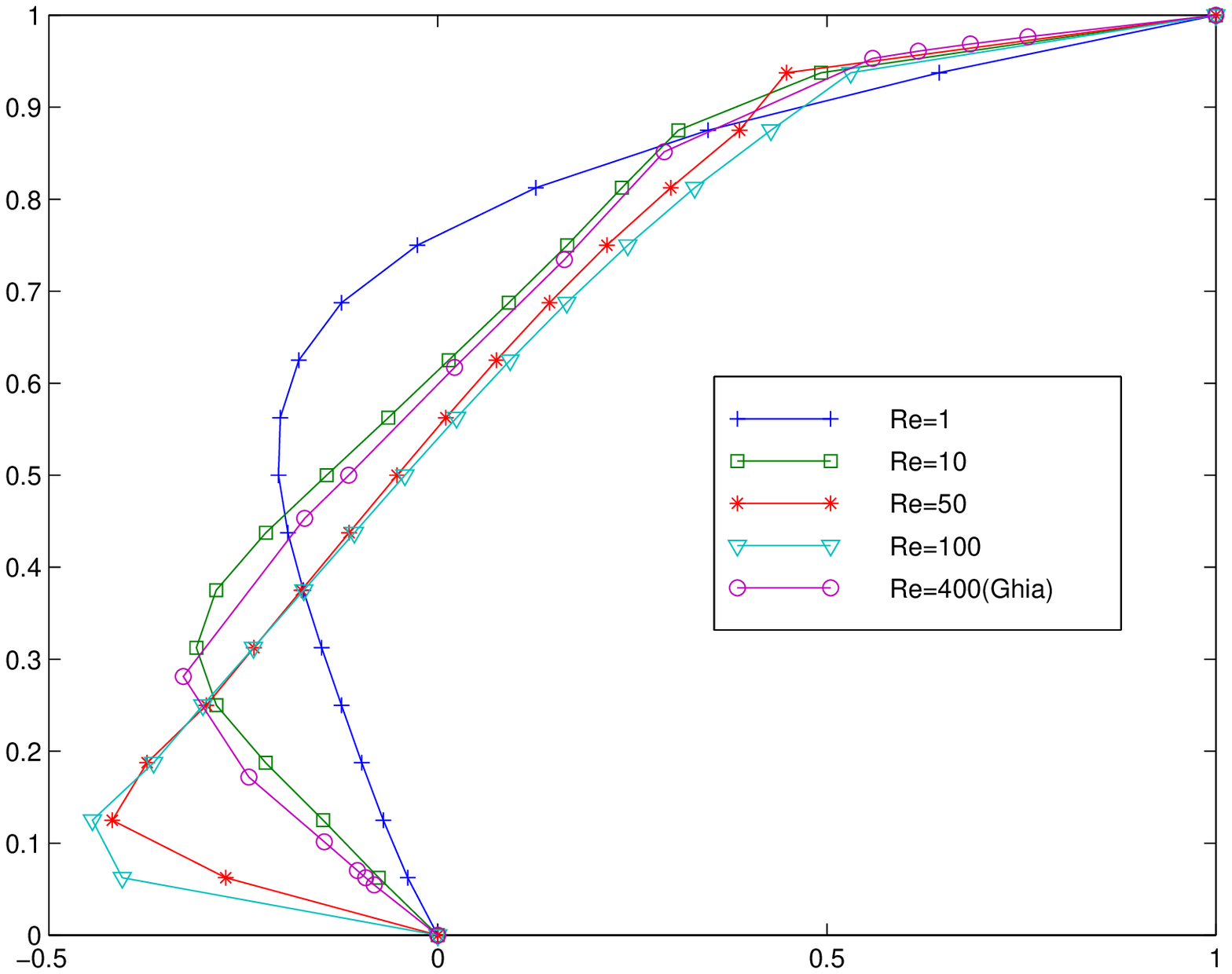,
	  width=2.5in,height=2.5in} 
	       } 
\centerline{
\epsfig{figure=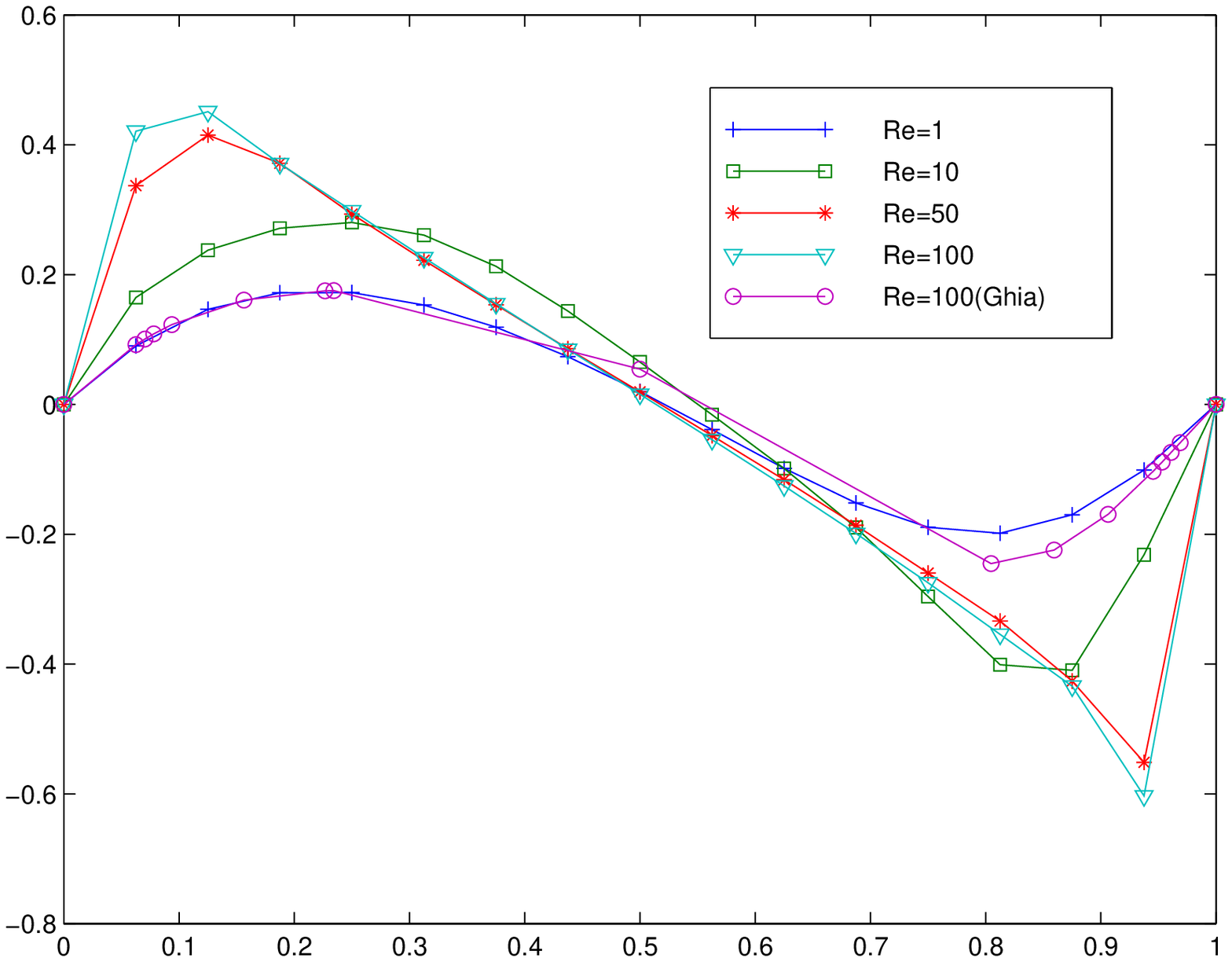,
	  width=2.5in,height=2.5in}
\epsfig{figure=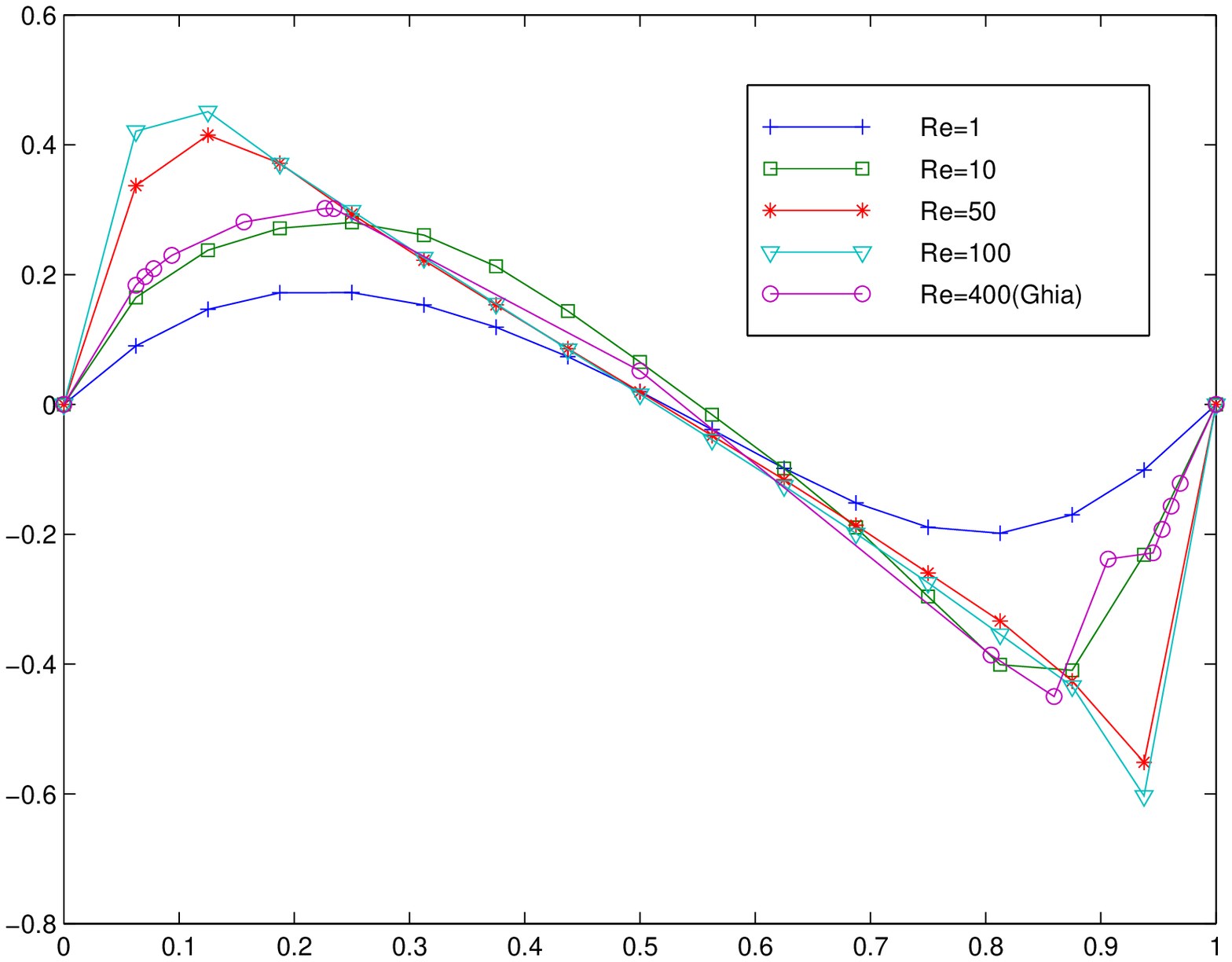,
	  width=2.5in,height=2.5in}	  
	  }  
  \caption{Cavity Problem :(above)$ u $-velocity lines through the
  vertical
  line $ x=0.5 $ for different Reynolds number and Ghia lines for
  Re=400, (below) 
  $ v $-velocity lines through the horizontal line $ y=0.5 $ for
  different 
  Reynolds number and Ghia lines for Re=400. (courtesy U.
  Ghia~\cite{ggs82}) }\label{fig4}
\end{figure}

We consider a domain $\Omega = [0, 1] \times [0, 1]$ with
no-slip boundary conditions, i.e., $u = v = 0$ in all
boundaries except $y=1$, where $u=1$. This problem has been
studied and addressed by many researchers including Ghia, Ghia, Shin
\cite{ggs82}, and J.E. Akin \cite{akin94}. The numerical 
computation of this example was obtained using a Sun Ultra 2 
with 2 200 MHz ultrasparc processor running Solaris 2.5.1. 
Bogner-Fox-Schmit elements are used with $17 \times
17$ grid points on the coarse mesh and $33 \times 33$ grid points on
the fine  mesh. The streamlines for Re = 1, 10, 50, 100 are plotted in
Figure~\ref{fig3}. Figure~\ref{fig4} shows the $u$-velocity lines
through the vertical line $x=0.5$ and $v$-velocity lines through 
the horizontal line $y = 0.5$.
% ----------------figure 1 --------------
%--------------Table 4 ---------------
\begin{longtable}[H]{|c|c|c|c|c|c|c|}
\multicolumn{7}{c}{ } \\ 
\hline  
$ Re $ & (H,h)  & noc  & nof & $\parallel
\psi-\psi^{h}\parallel_{0,h}$ & $\parallel
\psi-\psi^{h}\parallel_{1,h}$ & $\parallel
\psi-\psi^{h}\parallel_{2,h}$  \\
\hline \hline 
10 & $(\frac{1}{16},\frac{1}{32})$ & 61,66,43 & 576 &
4.21$\times10^{-5}$ & 4.84$\times10^{-3} $ & 5.27$\times10^{-2}$ \\
\hline
50 & $(\frac{1}{16},\frac{1}{32})$ & 54,83,173, & 1676 &
1.11$\times10^{-4}$ & 4.82$\times10^{-3}$ & 5.19$\times10^{-2}$ \\
 &  & 292,82 & & & &  \\
\hline
100 & $(\frac{1}{16},\frac{1}{32})$ & 53,254,421, & 4356 &
4.44$\times10^{-4}$ & 5.00$\times10^{-3}$ & 5.11$\times10^{-2}$ \\
 &  & 668,1287 & & & &  \\ 
\hline
200 & $(\frac{1}{16},\frac{1}{32})$ & 53,254,421, & 4356
& 3.06$\times10^{-4}$ & 4.52$\times10^{-3}$ & 5.07$\times10^{-2}$ \\
 &  & 668,1287,454 & & & &  \\
\hline
1000 & $(\frac{1}{16},\frac{1}{32})$ & 53,254,421, &
4356 & 9.62$\times10^{-4}$ & 1.08 $\times10^{-2}$ & 1.83
$\times10^{-1} $ \\ 
 &  & 668,1287,1156 & & & &  \\
\hline
2000 & $(\frac{1}{16},\frac{1}{32})$ &
53,254,421, & 4356 & 1.93  $\times10^{-3}$ &
2.59 $\times10^{-2}$ & 4.93 $\times10^{-1}$ \\
 &  & 668,1287,1156, & & & &  \\ 
&  & 1156,1156 & & & &  \\
\hline 
\caption{Two level for test problem.   where noc = number of 
Bicgstab iterations in coarse , nof = number of Bicgstab iteration
in fine}\label{table4}
\end{longtable}
%--------------Table 4 ---------------
% ----------------figure 1 --------------
% ------------ figure 2 --------------
% ----------------figure  --------------
% ----------------figure  3 & 4--------------
% -----------------------------------------------
%\newpage

\end{document}